\documentclass[12pt,a4paper,twoside]{amsart}
\usepackage{a4,amscd,amssymb,amsmath,amsthm,eucal,fancyheadings,xypic}
\usepackage[all]{xy}
\usepackage[english]{babel}
\numberwithin{equation}{section}
\newtheorem{thm}{Theorem}[section]

\newtheorem{prop}[thm]{Proposition}
\newtheorem{cor}[thm]{Corollary}

\newcommand{\z}{\zeta}
\newcommand{\al}{\alpha}
\newcommand{\be}{\beta}

\newcommand{\e}{\epsilon}
\newcommand{\la}{\lambda}
\newcommand{\w}{\omega}
\newcommand{\Z}{\mathbb Z}
\newcommand{\Q}{\mathbb Q}

\newcommand{\V}{\mathcal V}
\newcommand{\m}{\mbox{ mod }}
\newcommand{\M}{\mathcal M}

\DeclareMathOperator{\pic}{Pic}
\DeclareMathOperator{\cl}{Cl}
\DeclareMathOperator{\kernel}{ker}

\DeclareMathOperator{\Image}{Im}

\DeclareMathOperator{\cha}{Char}
\DeclareMathOperator{\Gal}{Gal}

\begin{document}
\title[The Kervaire-Murthy Conjectures]{Fine Structure of
Class Groups
$\cl^{(p)}\Q(\z_n)$
and the Kervaire--Murthy Conjectures II}
\author{Ola Helenius}
\address{Department of Mathematics, Chalmers University of
Technology
 and G\"{o}teborg University, SE-41296 G\"{o}teborg, Sweden}
\email{olahe@math.chalmers.se, astolin@math.chalmers.se}

\author{Alexander Stolin}
\address{Department of Mathematics, Chalmers University of
Technology and
G\"{o}teborg University, SE-41296 G\"{o}teborg, Sweden}

\date{}
%
%
\begin{abstract}
There is an Mayer-Vietoris exact sequence
\[
0\to V_n\to \pic \Z C_{p^n}\to \cl \Q (\z_{n-1})\times   \pic
\Z
C_{p^{n-1}}\to 0
\]
involving the Picard group of the integer group ring $\Z
C_{p^n}$
where $C_{p^n}$ is the cyclic group of order $p^n$ and
$\zeta_{n-1}$ is a primitive $p^n$-th root of unity. The
group
$V_n$ splits as $V_n\cong V_n^+\oplus V_n^-$ and $V_n^-$ is
explicitly known. $V_n^+$ is a quotient of an in some sense
simpler group $\mathcal{V}_n$. In 1977 Kervaire and Murthy
conjectured that for semi-regular primes $p$, $V_n^+ \cong
\mathcal{V}_n^+ \cong \cl^{(p)}(\Q (\zeta_{n-1}))\cong
(\mathbb{Z}/p^n
\mathbb{Z})^{r(p)}$, where $r(p)$ is the index of regularity
of
$p$. Under an extra condition on the prime $p$, Ullom
calculated
$V_n^+$ in 1978 in terms of the Iwasawa invariant $\lambda$
as
$V_n^+ \cong (\mathbb{Z}/p^n \mathbb{Z})^{r(p)}\oplus
(\mathbb{Z}/p^{n-1} \mathbb{Z})^{\lambda-r(p)}$.

In the previous paper we
proved that for all semi-regular primes, $\mathcal{V}_n^+
\cong
\cl^{(p)}(\Q (\zeta_{n-1}))$ and that these groups are
isomorphic to
\[(\mathbb{Z}/p^n \mathbb{Z})^{r_0}\oplus (\mathbb{Z}/p^{n-1}
\mathbb{Z})^{r_1-
r_0}
\oplus \hdots \oplus (\mathbb{Z}/p
\mathbb{Z})^{r_{n-1}-r_{n-2}}
\]
for a certain sequence $\{r_k\}$ (where $r_0=r(p)$).
Under Ulloms extra condition it was proved that \[V_n^+ \cong
\mathcal{V}_n^+ \cong \cl^{(p)}(\Q(\z_{n-1})) \cong
(\mathbb{Z}/p^n
\mathbb{Z})^{r(p)}\oplus (\mathbb{Z}/p^{n-1}
\mathbb{Z})^{\lambda-r(p)}.\] In the present paper
we prove that Ullom's extra
condition is valid for all
semi-regular primes and it is hence shown that the above
result
holds for all semi-regular primes.

\end{abstract}
\maketitle

The present paper is a continuation of the paper
\cite{article-1-02}
of the authors and we refer you there for a more thorough
introduction.

Let $p$ be an odd prime,  $C_{p^n}$ denote the cyclic group
of
order $p^n$ and let $\z_{n}$ be a primitive $p^{n+1}$-th root
of
unity. In this paper we work on problems related to $\pic (\Z
C_{p^n})$. Our methods also lead to the calculation of the
$p$-part of the ideal class group of $\Z[\z_n]$. Recall that
calculating Picard groups for a group ring like the one above
is
equivalent to calculating $K_0$ groups.

There is a well known exact sequence involving the Picard
group of
$\Z C_{p^n}$ which was for example presented by Kervaire and
Murthy in \cite{K-M}. The sequence, which is based on the
$(*,\pic)$-Mayer-Vietoris exact sequence associated to a
certain
pullback of rings, reads
\begin{equation}\label{ZC-MV}
0\to V_n \to \pic \Z C_{p^{n+1}} \to \pic \Z C_{p^{n}}\times
\cl \Q(\z_n) \to 0.
\end{equation}

In \cite{article-1-02} we observe that $\pic (\Z
C_{p^{n}})\cong
\pic A_n$, where
\[
A_n:=\frac{\Z[x]}{\big( \frac{x^{p^n}-1}{x-1} \big)}.
\]
We look at the pullback
\begin{equation}\label{pullbackdiagram}
\xymatrix@=40pt{ A_{n+1} \ar[r]^{i_{n+1}} \ar[d]_{j_{n+1}} &
\Z[\z_{n}] \ar[dl]^{N_{n}}
\ar[d]^{f_{n}}\\
A_{n} \ar[r]^{g_{n}}         & \frac{A_{n}}{pA_{n}}=:D_{n} }
\end{equation}
where the map $N_{n}$ is constructed so that the lower right
triangle of the diagram is commutative. From this we get an
exact
sequence equivalent to \ref{ZC-MV} where $V_n$ can be
represented
as
\[
V_n = \frac{D_n^*}{g_n(A_n^*)}.
\]
Here $R^*$ denotes the group of units in a ring $R$. Using
\ref{pullbackdiagram} and the map $N_n$ we can construct an
embedding of $\Z[\z_{n-1}]^*$ into $A_n^*$ which we by abuse
of
notation consider an identification. This allows us to define
\[
\V_n:=\frac{D_n^*}{g_n(\Z[\z_{n-1}]^*)}.
\]
There is an action of $G_n:=\Gal(\Q(\z_n)/\Q)$ on all
involved
groups so we can use complex conjugation $c$ and for each
involved
(multiplicative) group  $G$ define $G^+$ and $G^-$ as the
subgroups invariant ($c(g)=g$) respectively anti-invariant
($c(g)=g^{-1}$) under $c$. Since both $V_n$ and $\V_n$ are
$p$-groups and $c$ is even we get a splitting
$V_n=V_n^+\oplus
V_n^-$ and $\V_n=\V_n^+\oplus \V_n^-$. It turns out that
$\V_n^+$
is isomorphic to its counterpart in \cite{K-M} (also denoted
by
$\V_n^+$).

Recall that a prime $p$ is semi-regular when $p$ does not
divide
the order of the ideal class group of the maximal real
subfield
$\Q(\z_0+\z_0^{-1})$ of $\Q(\z_0)$. A non semi-regular prime
has
yet to be found and it is an old conjecture by Vandiver that
all
primes are semi-regular. We also recall that the index of
regularity $r(p)$ is defined as the number of Bernoulli
numbers
$B_2, B_4,\hdots,B_{p-3}$ with numerators (in reduced form)
divisible by $p$. Kervaire and Murthy conjecture in
\cite{K-M}
that for semi-regular primes:
\begin{eqnarray}
V_n^+ &=& \V_n^+ \label{eq:conj1}\\
\cha \V_n^+ &=& \cl^{(p)} \Q(\z_{n-1}) \label{eq:conj3}\\
\cha V_n^+ &\cong&
\big(\frac{\Z}{p^n\Z}\big)^{r(p)},\label{eq:conj2}
\end{eqnarray}
Results mainly from Iwasawa show that for big enough $n$,
\[
|\cl^{(p)}(\Q(\z_{n})|=p^{\la n + \nu}
\]
and the constants $\la$ and $\nu$ are called Iwasawa invariants.
Resulting from as splitting of $\cl^{(p)}(\Q(\z_{0}))$ with
respect to idempotents there are also an Iwasawa invariants $\la_i$
for each component $e_i S_n$ of $S_n=\cl^{(p)}(\Q(\z_{n})$. In 1978 Ullom
showed in \cite{U2} that if each $\la_{1-i}$ satisfy $
1\leq\la_{1-i}\leq p-1 $, then
\begin{equation}\label{eq:Ullom2}
V_n^+\cong (\frac{\Z}{p^n \Z})^{r(p)} \oplus
(\frac{\Z}{p^{n-1}\Z})^{\la-r(p)}.
\end{equation}

In our previous paper \cite{article-1-02} we prove a number
of
results regarding the Kervaire--Murthy conjectures. Before
presenting them we need a some definitions. For $n \geq 0$
and
$k\geq 0$, define
\[
U_{n,k}:=\{ real\,\,\e \in \Z[\z_n]^* : \e \equiv 1 \m
\la_n^{k}\},
\]
where $\la_n=(\z_n-1)$ is the prime above $(p)$ in
$\Z[\z_n]$. Let
$U^{p}$ denote the group of $p$-th powers of elements of the
group
$U$. Note that we in this paper sometimes use the notation
$R^n$
for $n$ copies of the ring (or group) $R$. The context will
make
it clear which one of these two things we mean. Similarly the
context should make it clear wether an indexed $\la$ means an
Iwasawa invariant or a prime ideal.

For $k=0,1,\hdots$, define $r_k$ by
\[
|U_{k,p^{k+1}-1} / (U_{k,p^{k}+1})^{p}|=p^{r_k}.
\]
It turns out that $r_0=r(p)$ (see for example
[B-S]). Our main results from
\cite{article-1-02}
are:
\begin{thm}\label{thm:article-1-02}
For semi regular primes,
\[
\cha \V_n^+ = \cl^{(p)} \Q(\z_{n-1})\cong
\big(\frac{\Z}{p^{n}\Z}\big)^{r_0}\oplus
\big(\frac{\Z}{p^{n-1}\Z}\big)^{r_1 - r_0}\oplus \hdots
\oplus
\big(\frac{\Z}{p\Z}\big)^{r_{n-1} - r_{n-2}}.
\] If Ulloms assumption holds, then $r_k=\la$ for all $k\geq
1$,
$\nu=r(p)=r(0)$ and
\begin{equation}\label{eq:main}
\cha V_n^+ \cong \cha
\V_n^+\cong \cl^{(p)} \Q(\z_{n-1}) \cong
\big(\frac{\Z}{p^{n}\Z}\big)^{r(p)}\oplus
\big(\frac{\Z}{p^{n-1}\Z}\big)^{\la-r(p)}.
\end{equation}
Moreover, $\la=r(p)$ is equivalent to that all three
Kervaire--Murthy conjectures hold and if $\la$ equals r(p),
then
$\nu$ equals $r(p)$.
\end{thm}

In the present paper we prove that Ulloms assumption on the
Iwasawa invariants $\la_{1-i}$ above is true for all semi-regular
primes. Before we can explain exactly what we prove we need some
more notation. We remind the reader that
$S_n$ denotes the $p$-part of the ideal class
group of $\Q(\z_n)$. One can find a set of mutually orthagonal
idempotents, $\e_i$, such that $S_n=\bigoplus \e_i S_n$. In 
\cite{W} it was proved that 
\[
\e_i S_n\cong\frac{\Z_p[[T]]}{((1+T)^{p^n}-1,f_i(T))},
\]
where $f_i(T)=a_0+a_1 T+a_2 T^2+\hdots$ is a power series
satisfying $f_i ((1+p)^s -1)=L_p (s,\omega^{1-i})$. Here
$L_p (s,\chi)$ is the $p$-adic $L$-function with a Dirichlet
character $\chi$ defined for instance in \cite{W} and 
$\omega (a)$ is a $p$-adic Dirichlet character of conductor
$p-1$. 

For the
constant term we have $a_0=-B_{1,\w^{-i}}=L_p (0,\w^{1-i})$, where
$B_{1,\w^i}$ is a generalized Bernoulli number
(again, see \cite{W}).

The Iwasawa invariants $\la_i$ turn out to be the first exponent
such that $p\not{|}a_{\la_i}$. Let ${\mathcal O}$ be a finite extension of
$\Z_p$ and ${\M}$ its maximal ideal. A polynomial $h\in
{\mathcal O}[T]$ is called distinguished if it has leading
coefficient $1$ and all other coefficients belong to ${\M}$. It
is known (see for instance Proposition 7.2 in \cite{W})
that for a distinguished polynomial $h$,
\[
\frac{{\mathcal O}[[T]]}{(h(T))}\cong \frac{{\mathcal O}[T]}{(h(T))}.
\]
In our case, using Weierstrass preparation theorem one can find a
distinguished polynomial
\[
g_i(T)=a_0'+a_1'T+ \hdots
+a_{\la_i-1}'T^{\la_{i-1}}+T^{\la_{i}}
\]
and an invertible series $u_i(T)$ such that $f_i(T)=g_i(T)u_i(T)$.
This representation is unique. Using this we get
\[
\e_i S_n\cong \frac{\Z_p[[T]]}{((1+T)^{p^n}-1,g_i(T))}\cong
\frac{\Z_p[T]}{((1+T)^{p^n}-1,g_i(T))}.
\]
Recall that we are interested in evaluating $\la_i$. First, for
$n=0$ we get that
\[
\e_i S_0 \cong \frac{\Z_p[T]}{(T,g_i(T))}\cong \frac{\Z_p}{(a_0)}
\cong \frac{\Z_p}{(a_0')}.
\]
From our previous results on $S_n$ we know that $\e_i S_n\cong
\Z/p\Z$, so $a_0'=pu$ for some unit $u\in \Z_p $. Hence $g_i$ is an
Eisenstein polynomial and hence irreducible. Now consider the case
$n=1$. We get
\[
\e_i S_1\cong \frac{\Z_p[T]}{((1+T)^{p}-1,g_i(T))}.
\]
Choose $\be_i$ such that $g_i(\be_i)=0$. Then,
\[
\e_i S_1\cong \frac{\Z_p[\be_i]}{((1+\be_i)^{p}-1)}.
\]
Suppose $\la_i\geq p$. The field $\Q_p (\be_i)$ completely
remifies over $\Q_p$ and has degree $\la_i$.
Therefore $(\be_i)^{\la_i}=(p)$  and we 
see that $(1+\be_i)^{p}-1=u\be_i^p$ for some unit $u\in\Z_p [\be_i]$.
Then we get that
\[
\frac{\Z_p[\be_i]}{((1+\be_i)^{p}-1)}=\frac{\Z_p[\be_i]}{(\be_i^p)}
\] 
and multiplication by $p$ annihilates this factor-ring.
Therefore for some $k$ we have: 
\[
\frac{\Z_p[\be_i]}{((1+\be_i)^{p}-1)}\cong (\Z/p\Z)^k.
\]
So, if we deduce from our previous results that there are
elements of order $p^2$ in $\e_i S_1$, it will contradict
to the assumption that $\la_i\geq p$.
We will hence prove the following two theorems.
\begin{thm}
Let $p$ be a semi-regular prime and let $g_i$ be the distinguished
polynomial defined above, with $a_0'$ being the constant
coefficient. Then we have
 \begin{enumerate}
  \item $p^2\not{|}\ a_0'$
  \item $g_i(T)$ is an Eisenstein
  polynomial of degree strictly less
  than $p$.
 \end{enumerate}

\end{thm}

\begin{thm}For semi-regular primes,
\begin{enumerate}
  \item $\la_{1-i}$ satisfy $1\leq\la_{1-i}\leq p-1$.
  \item $r_k=\la$ for all $k=1,2,3\hdots$.
  \item $V_n^+\cong \V_n^+\cong \cha S_{n-1}
        \cong \big(\frac{\Z}{p^n \Z}\big)^{r(p)}
        \oplus \big(\frac{\Z}{p^{n-1} \Z}\big)^{\la -r(p)}$
\end{enumerate}
\end{thm}
{\bf Remark.} The above yields for semi-regular $p$ that
$\nu=r(p)$.
\vskip0.3cm
Let us recall that 
\[
S_1\cong (\Z/p^2\Z)^{r(p)}\oplus (\Z/p\Z)^{r_1-r(p)}
\]
The usual norm map induces an epimorphism $N:S_1\to S_0$.
\begin{prop}
\[
\kernel N\cong (\Z/p\Z)^{r_1}
\]
\end{prop}
\begin{proof}
Let $A$ be any finite abelian group and 
let us denote by $\cha A$ or $A^{\times}$ the group of its characters.
Clearly any homomorphism $f: A\to B$ induces a dual homomorphism
$f^* :B^{\times}\to A^{\times}$.

In the proof of Theorem 2.14 in [H-S] we constructed an embedding
$\al_1:\V_0^+\cong S_0^{\times}\to S_1^{\times}\cong \V_1^+$. The map $\al_1$
is induced by the canonical embedding $\Q(\z_0)\to \Q(\z_1)$
and then clearly $N^*=\al_1$. Then we get that
\[
\kernel N\cong\cha \frac{\V_1^+}{\Image (\al_1)}
\]
Therefore we have to prove that 
\[
\frac{\V_1^+}{\Image (\al_1)}\cong (\Z/p\Z)^{r_1}
\]
For this we recall that we also have a surjection 
$\pi_1:\V_1^+\to\V_0^+$ [Proposition 2.12, H-S]. Moreover,
it was proved in the proof of Theorem 2.14 in [H-S] that
$\al_1 (\pi_1 (a))=a^p$ for any $a\in \V_1^+$.
The latter implies that 
\[
\frac{\V_1^+}{\Image (\al_1)}=
\frac{\V_1^+}{(\V_1^+)^p}\cong (\Z/p\Z)^{r_1}
\]
\end{proof}
\begin{cor}
$\e_i S_1$ contains elements of order $p^2$.
\end{cor}
\begin{proof}
It is known that $N$ maps $\e_i S_1$ onto $\e_i S_0$
(see for instance [W]). Since $S_1$ has $r_1$ generators
and $\kernel N\cong (\Z/p\Z)^{r_1}$ it follows that
any preimage of non-zero $a\in S_0$ has order $p^2$ and
hence, $\e_i S_1$ contains an element of order $p^2$.
This completes the proofs of Corrolary 0.5 and Theorems 0.2, 0.3.
\end{proof}
\vskip0.3cm
{\bf Final Remark.} The following result was proved in [W].
\begin{thm}
Suppose $p$ is semi-regular. Let an even index $i$ be such that
$2\leq i\leq p-3$ and $p|B_i$ ($B_i$ is the corresponding 
ordinary Bernoulli
number). If
\[
B_{1,\w^{i-1}}\not{\equiv} 0\m p^2
\]
and
\[
\frac{B_i}{i}\not{\equiv}\frac{B_{i+p-1}}{i+p-1}\m p^2
\]
then for all $n\geq 0$
\[
S_n\cong (\Z/p^{n+1}\Z)^{r(p)}
\]
For semi-regular $p$ the above yields
\[
\la=\nu=r(p)
\]
\end{thm}
It was written in a remark after the result that the above
incongruences hold for all $p<4000000$ but there does not seem
to be any reason to believe this in general.

Our results show that the first incongruence is valid for all
semi-regular primes as well as $\nu=r(p)$. 
So we may hope that the second incongruence above obtained numerically
also is valid in some generality.

\end{document}